\documentclass[a4paper]{amsart}
\usepackage[utf8]{inputenc}

\usepackage{amsmath,amsfonts,amssymb}
\usepackage{amsthm}
\usepackage{latexsym, xspace, enumerate}
\usepackage[mathscr]{eucal}
\usepackage[all]{xy}
\usepackage{hyperref}
\usepackage{amscd}
\usepackage{mathtools}
\usepackage{vmargin}
\usepackage{rotating}
\usepackage{tikz,comment}
\usepackage{comment}

\usepackage{stackengine,scalerel}
\stackMath

\newtheorem{theorem}{Theorem}[section]
\newtheorem{lemma}[theorem]{Lemma}

\newtheorem{thmx}{Theorem}
\newtheorem{corx}[thmx]{Corollary}

\theoremstyle{definition}

\def\BS{\mathrm{BS}}

\makeatletter
\def\Ddots{\mathinner{\mkern1mu\raise\p@
\vbox{\kern7\p@\hbox{.}}\mkern2mu
\raise4\p@\hbox{.}\mkern2mu\raise7\p@\hbox{.}\mkern1mu}}

\title{A new family of sofic one-relator groups}
\author{Federico Berlai}

\address[Federico Berlai]{Department of Mathematics, University of the Basque Country UPV/EHU, Barrio Sarriena s/n, 48940 Leioa, Spain}
\email[Federico Berlai]{federico.berlai@ehu.eus}

\begin{document}

\begin{abstract}
    We provide an infinite family of sofic one-relator groups that are not residually finite nor residually solvable. The proof is essentially different from the one in \cite{Bannon} as it does not use Magnus' decompositions.
\end{abstract}

\maketitle

Throughout the paper we use the notation $x^y:=y^{-1}xy$ and $[x,y]:=x^{-1}y^{-1}xy$.
%Recall the following open problems:
%\begin{openproblem}
%\begin{enumerate}
%\item %\cite[Problem 18.6]{Kou} 
%\emph{Are all one-relator groups residually amenable?}
%\item %\cite[Open Question 4.9]{Pestov} 
%\emph{Are all one-relator groups sofic?}
%\end{enumerate}
%\end{openproblem}
%The first question is due to Goulnara Arzhantseva \cite[Problem 18.6]{Kou}, while the second is a variation of a question of Nate Brown \cite[Open Question 4.9]{Pestov}. 

%Recall that a group is residually amenable if amenable quotients separate elements, while a group is sofic if it is a subgroup of a metric ultraproduct of finite symmetric groups equipped with the normalised Hamming distance \cite{Pestov}. From the given definitions a connection between the two concepts seems obscure; nevertheless, it is known that residual amenability implies soficity. Thus, a positive answer to the first problem implies a positive answer to the second. Needless to say, it is a major open problem to see if all groups are sofic.

%\smallskip
In this short note, we focus on a family of one-relator groups introduced by Baumslag, Miller~III and Troeger~\cite{BM3T} in 2007. These one-relator groups generalize a classical example of Baumslag~\cite{Baum69} from 1969, and they are not residually finite nor residually solvable. To introduce this family of one-relator groups, let $r,w\in F(a,b)$ be two elements that do not commute (more generally, one could consider a free group on more than two generators, but for the sake of simplicity we will focus on the two-generated case here). Define
\begin{equation}\label{generalized_Baumslag_eq}
G_{r,w}:= \langle a,b \mid r=[r,r^w]\rangle= \langle a,b \mid r^2=r^{r^w}\rangle,
\end{equation}
In \cite{BM3T} it is proved that $G_{r,w}$ has the same finite quotients as the group $\langle a,b \mid r\rangle$, and that $G_{r,w}$ is not residually finite (because, for instance, the non-trivial element $r$ of $G_{r,w}$ has trivial image in all finite quotients).
Such groups are not residually solvable either, as the element $r$ belongs to all derived subgroups of $G_{r,w}$.
In the case $r=a$ and $w=b$ one recovers the classical Baumslag group
\begin{equation}\label{eq_BG_presentation}
    B:=\langle a,b\mid a=[a,a^b]\rangle=\langle a,b\mid a^2=(b^{-1}ab)^{-1}a(b^{-1}ab)]\rangle.
\end{equation}
This group was considered by Gilbert Baumslag in 1969 to provide (yet another) example of a non residually finite one-relator group \cite{Baum69}. 
Bannon \cite{Bannon} (compare also \cite{BannonNoblett}) noticed that, despite $B$ not being residually finite nor residually solvable, it is a sofic group. Indeed, looking at Magnus' decomposition for such one-relator group one sees that 
$B$ is an HNN-extension of a Baumslag-Solitar group $\BS(1,2)$ amalgamating a cyclic subgroup, that is, $B=\BS(1,2)\ast_\mathbb{Z}$. 
%\[
%\begin{split}
%B&=\langle a,b\mid a^2=(b^{-1}ab)^{-1}a(b^{-1}ab)]\rangle=
%\langle a,t,b\mid a^2=a^t,t=b^{-1}ab\rangle=\BS(1,2)\ast_\mathbb{Z}.
%\end{split}
%\]
%is an HNN extension with metabelian (hence sofic) base group and amenable amalgamated subgroup. 
Such HNN extensions are known to be sofic \cite{ElSz, Pau} because the base group is sofic and the amalgamated subgroup is amenable. This was the first example of a sofic one-relator group that is not residually finite nor residually solvable, with more examples considered in \cite{BannonNoblett}, applying the same strategy.

As noticed in \cite{BannonNoblett}, soficity of the groups $G_{a,b^n}$ for $n>1$ does not follow from the argument of \cite{Bannon,BannonNoblett},
because their Magnus' decomposition is
$(\BS(1,2)\ast F_{n -1})\ast_{F_{ n}}$,
where (for $n>1$) a non-amenable free group appears as amalgamated subgroup. For such HNN extensions there is no general result implying soficity. The definition of sofic group will not play a role here, hence we refer to \cite{CecSil,Pestov} for excellent introductions to these groups. We will just underline that it is an open problem to see if all groups are sofic \cite[Open Question 3.8]{Pestov}, or if all one-relator groups are sofic \cite[Open Question 4.9]{Pestov}.

\smallskip
In what follows, we will mostly focus on the case $r=a$ and plug in exponents $l,k\in \mathbb{Z}$ into the relator in a fashion similar to Baumslag-Solitar groups, so that the relation reads $r^k=(r^l)^{r^w}$, to cover a slighly more general family of groups (one could go even further and consider relators as done in~\cite{Bru}, but we will not do this here). Given two elements $r,w\in F(a,b)$ that do not commute and two non-zero integers $l$ and $k$, we define
\begin{equation}\label{generalized2}
G_{r,w}(l,k):= \langle a,b,\ldots \mid (r^l)^{r^w}=r^k\rangle.
\end{equation}
As justified in the following section, many of the one-relator groups from Equation \eqref{generalized2} are not residually finite nor residually solvable.
Whenever $l=1$ and $k=2$, we simplify the notation and write directly $G_{r,w}$ instead of $G_{r,w}(1,2)$ (these are exactly the groups recalled in Equation~\eqref{generalized_Baumslag_eq}).~We~prove:
\begin{thmx}
For all non-zero integers $n,l,k\in \mathbb{Z}$, the group $G_{a,b^n}(l,k)$ is sofic.% If $\lvert l\rvert \neq \lvert k\rvert$ then it is not residually finite, and if $k=l+ d$ where $d$ is an integer divisor of $l$ then it is not residually solvable.}
\end{thmx}
This settles the open problem raised in \cite{BannonNoblett} (see \cite[Theorem 2]{BannonNoblett} and the discussion before that).
A careful choice of the integers $l$ and $k$ immediately implies:
\begin{corx}
For all $n\in \mathbb{Z}\setminus\{0\}$ and $l\in \mathbb{Z}\setminus\{0,-1\}$ the group $G_{a,b^n}(l,l+1)$ is sofic but not residually finite nor residually solvable.
\end{corx}

Theorem A can be used iteratively to generate more sofic one-relator groups that are not residually finite nor residually solvable:
\begin{thmx}
For all non-zero integers $n,l,k\in \mathbb{Z}$ the one-relator group
$G_{a,b^{-n}ab^n}(l,k)$ %which also have the presentation $\langle a,b,t\mid a^2=a^{a^t}, t=b^{-n}ab^n\rangle$,
is sofic.
\end{thmx}
As in Corollary B, correctly choosing the integers $l$ and $k$ produces groups that are sofic but not residually finite nor residually solvable.

\medskip
It is not evident if all groups among the one-relator groups introduced by Baumslag, Miller III and Troeger can be shown to be sofic using Magnus' decompositions and the argument used in this paper (or by other methods). For instance, is the group $G_{a,b^{-1}ab^2}$ sofic?

\subsubsection*{Acknowledgements}
The author is supported by the Spanish Government, grant PID2020-117281GB-I00, partly by the European Regional Development Fund (ERDF), and the Basque Government, grant IT1483-22. He is indebted to Goulnara Arzhantseva for suggesting to work on this project during his PhD studies, back in 2013.
He is very thankful to Marco Linton for insightful conversations on this topic and for pointing out a gap in a proof in a previous version of this note, and to an anonymous referee for suggesting a slightly shorter/cleaner proof of Theorem A.

\section{Preliminary lemmas}

%In general, if $l=1$ and $k=2$ then the group $G_{a,b^n}(1,2)$ is not residually solvable: the element $a=[a,b^{-n}ab^n]$ is non-trivial and belongs to all derived subgroups of $G_{a,b^n}(1,2)$. 

In this preliminary section we collect some observations that will be useful later, justifying our generalisation of the family of one-relator groups introduced by Baumslag, Miller III and Troeger~\cite{BM3T}. In what follows, $r$ and $w$ are elements of the free group $F(a,b)$ that do not commute.

We start by noticing that many of these groups are not residually finite. Indeed, if
$l=\pm1$ and $k=\pm2$ (or viceversa), 
then $G_{r,w}(l,k)$ is not residually finite. In particular, if
$(l,k)\in\{(1,2),(-1,-2),(2,1),(-2,-1)\}$, then the claim follows directly applying \cite[Theorem 1]{BM3T}. If signs are mixed, that is if $(l,k)\in\{(1,-2),(-1,2),(2,-1),(-2,1)\}$, the claim essentially follows from \cite[Theorem 1]{BM3T}, again, because \cite[Lemma 1]{BM3T} holds in this case: that is, if $x,y\in G$ are two elements of the same finite order and satisfy the relation $x^y=x^{-2}$, then the subgroup $\langle x,y\rangle$ of $G$ is trivial. 
If $l$ and $k$ are different in absolute value and both different from $\pm1$, then
the group $G_{r,w}(l,k)$ is not residually finite by \cite[Theorem B]{Meskin}.
These remarks and Lemma \ref{lemmanotRS}, joint with Theorem A, already imply Corollary B.

The groups $G_{r,w}$ considered in \cite{BM3T} are always non-residually finite, for any pair of non-commuting elements $r,w$ in a free groups. A full characterization of (non-) residual finiteness of the groups $G_{r,w}(l,m)$ is not given here, but it seems plausible to suppose that these groups also are never residually finite.

Many of these groups are also not residually solvable:
\begin{lemma}\label{lemmanotRS}
%If $k=l +d$ where $d$ is an integer dividing $l$, then the group $G_{r,w}(l,k)$ is not residually solvable.
If $k\neq 0\neq l$ and $k= l\pm1$, then the group $G_{r,w}(l,k)$ is not residually solvable.
\begin{proof}
Suppose that $k=l+1$. The defining relation $r^k=(r^l)^{r^w}$ can be rewritten as ${r=r^{-l}(r^l)^{r^w}}$, that is $r=[r^l,(r)^w]$ is an element of the derived subgroup. Consequently the elements $r^w$ and $r^l$ also lie in the derived subgroup. Inductively, $r$ is an element of the intersection of all derived subgroups of $G_{r,w}(l,k)$. This element is not trivial by the same argument of \cite[Lemma 1]{BM3T}; thus, $G_{r,w}(l,k)$ is not residually solvable.

If $k=l-1$, then the same argument shows that $r^{-1}$ is a non-trivial element in the intersection of all derived subgroups of $G_{r,w}(l,k)$.
%Suppose that $l=dl_1$ for some integer $l_1$. The defining relation $r^k=(r^l)^{r^w}$ can be rewritten as ${r^d=r^{-l}(r^l)^{r^w}}$, that is $r^d=[r^l,(r^l)^w]$ is an element of the derived subgroup. Thus $r^l=(r^d)^{l_1}$ also is an element of the derived subgroup, and consequently the conjugate $(r^l)^w$ too. Inductively, $r^d$ is an element of the intersection of all derived subgroups of $G_{r,w}(l,k)$. This element is not trivial by the same argument of \cite[Lemma 1]{BM3T}; thus, $G_{r,w}(l,k)$ is not residually solvable.
\end{proof}
\end{lemma}

Even if it will be of no use in this note, we record the following observations:

\begin{lemma}
Let $r,w\in F(a,b)$ be two non-commuting elements and $l,k\in \mathbb{Z}$. Then $G_{r,w}(l,k)\cong G_{r,w}(-l,-k)$ and $G_{r,w}(l,k)\cong G_{r^{-1},w}(k,l)$.
\begin{proof}
For the first isomorphism, notice that the relation $(r^l)^{r^w}=r^k$ implies that $\left((r^l)^{r^w}\right)^{-1}=\left(r^k\right)^{-1}$, that is 
$(r^{-l})^{r^w}=r^{-k}$. Therefore
\[\begin{matrix}
G_{r,w}(l,k)&\stackrel{\eta}{\longrightarrow} & G_{r,w}(-l,-k)& \stackrel{\xi}{\longrightarrow} & G_{r,w}(l,k)\\
a&\longmapsto &a & \longmapsto &a\\
b&\longmapsto &b & \longmapsto &b
\end{matrix}\]
are group homomorphisms such that $\xi\circ \eta$ and $\eta\circ \xi$ are both identity homomorphismsm. Therefore they are isomorphisms.

For the second isomorphism, notice that $(r^l)^{r^w}=r^k$ immediately implies that $r^l=(r^{k})^{(r^{-1})^w}$, and the claim follows taking inverses.
\end{proof}
\end{lemma}

\begin{lemma}
Let $w\in F(a,b)$ be a non-trivial element whose normal form starts and ends with non-trivial powers of $b$. Then $G_{a,a^nw}(l,k)\cong G_{a,w}(l,k)\cong G_{a,wa^n}(l,k)$.
\begin{proof}
This is easy to see by explicitly writing the relations of $G_{a,a^nw}(l,k)$ and of $G_{a,wa^n}(l,k)$ and reducing them to the relation of $G_{a,w}(l,k)$.
\end{proof}
\end{lemma}

\section{Proof of Theorem A}
We may assume without loss of generality that $n$ is positive, because the map defined on the generators as $a\mapsto a$ and $b\mapsto b^{-1}$ is an automorphism.
For the case $n=1$, note that $G:=G_{a,b}(l,k)$ is the HNN extension
\[\begin{split}
G_{a,b}(l,k)&=\langle a,b\mid (a^b)^{-1}a^la^b=a^k\rangle =
\langle a_0,a_1,b\mid a_1^{-1}a_0^la_1=a_0^k, b^{-1}a_0b=a_1\rangle =\BS(l,k)\ast_{\mathbb{Z}}.
\end{split}\]
The Baumslag-Solitar group $\BS(l,k)$ is residually solvable \cite{Krop90}, hence sofic. Thus, $G_{a,b}(l,k)$ is sofic \cite{ElSz, Pau}. 
Consider the the relator of the group of Equation \eqref{generalized2}.
As the generator $b$ has exponent-sum equal to zero in this relator, the group $G$ admits a surjective homomorphism onto $\mathbb{Z}=\langle z\rangle$ defined by the assignments $a\mapsto e$ and $b\mapsto z$. The kernel $K$ of this homomorphism is the normal closure in $G$ of the element $a$, that is
\begin{equation}\label{eqK}
\begin{split}
K&=\langle\langle a\rangle\rangle_G=\langle \{ a_i\}_{i\in \mathbb{Z}}\rangle
=\bigl\langle \{ a_i\}_{i\in \mathbb{Z}}\mid \{a_{i+1}^{-1}a_i^la_{i+1}=a_i^k\}_{i\in \mathbb{Z}}\bigr\rangle, \qquad \qquad a_i:=b^{-i}ab^{i}\quad\forall\, i\in \mathbb{Z}.
%&=\ldots \underset{a_0\mapsto a_0}{\ast}\langle a_0,a_1\mid a_{1}^{-1}a_0^la_{1}=a_0^k\rangle\underset{a_1\mapsto a_1}{\ast} \langle a_1,a_2\mid a_{2}^{-1}a_1^la_{2}=a_1^k\rangle \underset{a_2\mapsto a_2}{\ast}\ldots
\end{split}
\end{equation}
This kernel is sofic, being a subgroup of the sofic group $G_{a,b}(l,k)$.

For the general case, consider $n\geqslant 2$ and the group $G_{a,b^n}(l,k)$. Also in this case the exponent-sum of the generator $b$ in the relator is equal to zero. Therefore $G_{a,b^n}(l,k)$  admits a surjective homomorphism onto $\mathbb{Z}$ with kernel the normal closure of $a$ in $G_{a,b^n}(l,k)$:
\[
K_n:=\langle\langle a\rangle\rangle= \langle \{ a_i\}_{i\in \mathbb{Z}}\mid \{a_{i+n}^{-1}a_i^la_{i+n}=a_i^k\}_{i\in \mathbb{Z}}\rangle.
\]
Grouping the generators modulo $n$, we see that $K_n$ is the free product of $n$ copies of the group $K$ appearing in Equation \eqref{eqK}. Therefore $K_n$ is sofic as well \cite[Theorem 1]{ElSz-1}, and thus $G_{a,b^n}(l,k)$ too, being sofic-by-amenable \cite[Proposition~7.5.14]{CecSil}. Thus, the proof of Theorem A is completed.

\section{Proof of Theorem C}
We now deduce Theorem C building on the result and arguments of the previous section.
Let 
\[G:=G_{a,b^{-n}ab^n}(l,k)=\langle a,b\mid \bigl((b^{-n}ab^n)^{-1}a(b^{-n}ab^n)\bigr)^{-1}a^l\bigl((b^{-n}ab^n)^{-1}a(b^{-n}ab^n)\bigr)=a^k\rangle.\]
As the generator $b$ has exponent-sum equal to zero in the relator, $G$ admits a surjective homomorphism onto $\mathbb{Z}=\langle z\rangle$ defined on the generators by $a\mapsto e$ and $b\mapsto z$. The kernel of this homomorphism is the normal closure of $a$ in $G$: it is generated by the conjugates $a_i:=b^{-i}ab^i$, where $i\in \mathbb{Z}$, subject to the relations $(a_{n+i}^{-1}a_ia_{n+i})^{-1}a_i^l(a_{n+i}^{-1}a_ia_{n+i})=a_i^k$. Regrouping the generators modulo $n$, this kernel is isomorphic to the free product of $n$ copies of
\[\begin{split}
H:&=\langle \{g_i\}_{i\in \mathbb{Z}}\mid  
\{(g_{i+1}^{-1}g_i  g_{i+1})^{-1}g_i^l(g_{i+1}^{-1}g_i g_{i+1})=g_i^k\}_{i\in \mathbb{Z}}\rangle\\
&=\ldots \underset{g_0\mapsto g_0}{\ast}\langle g_0,g_1\mid (g_0^l)^{g_0^{g_1}}=g_0^k\rangle\underset{g_1\mapsto g_1}{\ast}
\langle g_1,g_2\mid (g_1^l)^{g_1^{g_2}}=g_1^k\rangle
\underset{g_2\mapsto g_2}{\ast}\ldots
\end{split}
\]
The group $G$ is locally sofic, and hence sofic \cite[Proposition 7.5.5]{CecSil}, being an iterated amalgamated product of $G_{a,b^{n}}(l,k)$ (which is sofic by Theorem A) amalgamated along cyclic (hence amenable) subgroups \cite{ElSz, Pau}.
Thus $K=H\ast \ldots \ast H$ is sofic as well \cite[Theorem 1]{ElSz-1}, and $G$ too, being an extension of the sofic group $K$ by the amenable quotient $\mathbb{Z}$.

%\subsubsection*{Further comments}
%A stronger result can be deduced applying verbatim the proof of \cite[Theorem 1]{BM3T}: the elementary amenable quotients of $G_{r,w}$ are exactly the elementary amenable quotients of $\langle a,b,\ldots\mid r\rangle$.

%It is not evident if all groups among the one-relator groups of Baumslag, Miller III and Troeger can be shown to be sofic using Magnus' decompositions and the argument used in this paper. In particular, is the group $G_{a,b^nab^nab^n}$ sofic?


\begin{thebibliography}{10}

\bibitem{Bannon} J. P. Bannon, \emph{A non-residually solvable hyperlinear one-relator group}. Proc. Amer. Math. Soc. 139 (2011), no. 4, 1409~-~1410;

\bibitem{BannonNoblett} J. P. Bannon, N. Noblett, \emph{A note on nonresidually solvable hyperlinear one-relator groups}. Involve 3 (2010), no. 3, 345~-~347;

\bibitem{Baum69} G. Baumslag, \emph{A non-cyclic one-relator group all of whose finite quotients are cyclic}. J. Austral. Math. Soc. 10 (1969), 497~-~498;


\bibitem{BM3T} G. Baumslag, C. F. Miller III, D. Troeger, \emph{Reflections on the residual finiteness of one-relator groups}. Groups Geom. Dyn. 1 (2007), no. 3, 209~-~219;

\bibitem{Bru} A. Brunner, \emph{On a class of one-relator groups}. Canadian J. Math. 32 (1980), no. 2, 414~-~420;

\bibitem{CecSil} T. Ceccherini-Silberstein, M. Coornaert, \emph{Cellular automata and groups}. Second edition. Springer Monographs in Mathematics. Springer, Cham, (2023);

%\bibitem{Chou} C. Chou, \emph{Elementary amenable groups}. Illinois J. Math. 24 (1980), no. 3, 396~-~407;

\bibitem{ElSz-1} G. Elek, E. Szab\'o, \emph{On sofic groups}. J. Group Theory 9, (2006), 161~-~171;

\bibitem{ElSz} G. Elek, E. Szab\'o, \emph{Sofic representations of amenable groups}. Proc. Amer. Math. Soc. 139 (2011), no. 12, 4285~-~4291;

%\bibitem{Kou} Kourovka notebook. \url{https://kourovka-notebook.org/wp-content/uploads/2024/06/20tkt.pdf}

\bibitem{Krop90} P. Kropholler, \emph{Baumslag-Solitar groups and some other groups of cohomological dimension
two}. Comment. Math. Helvetici 65 (1990), 547~-~558;

\bibitem{Meskin} S. Meskin, \emph{Nonresidually finite one-relator groups}. Trans. Amer. Math. Soc. 164 (1972), 105~-~114;

\bibitem{Pau} L. P\u aunescu, \emph{On sofic actions and equivalence relations}. J. Funct. Anal. 261 (2011), no. 9, 2461~-~2485;

\bibitem{Pestov} V. Pestov, \emph{Hyperlinear and sofic groups: a brief guide}.  Bull. Symbolic Logic 14 (2008), no. 4, 449~-~480.


\end{thebibliography}
\end{document}